# THE STRUCTURE OF A SOLVMANIFOLD'S HEEGAARD SPLITTINGS

DARYL COOPER AND MARTIN SCHARLEMANN

ABSTRACT. We classify isotopy classes of irreducible Heegaard splittings of solvmanifolds. If the monodromy of the solvmanifold can be expressed as

$$\begin{pmatrix} \pm m & -1 \\ 1 & 0 \end{pmatrix},$$

for some $m \geq 3$ (as always is true when the trace of the monodromy is $\pm 3$), then any irreducible splitting is strongly irreducible and of genus two. If $m \geq 4$ any two such splittings are isotopic. If $m = 3$ then, up to isotopy, there are exactly two irreducible splittings, their associated hyperelliptic involutions commute, and their product is the central involution of the solvmanifold.

If the monodromy cannot be expressed in the form above then the splitting is weakly reducible, of genus three and unique up to isotopy.

## 1. INTRODUCTION

The study of Heegaard splittings of 3-manifolds is now nearly a century old [He]. (See also [Prz] for a translation of the relevant parts). Such a splitting is deceptively simple to describe: a closed 3-manifold is regarded as the union of two handlebodies glued together along their boundaries. Although any 3-manifold can be described this way, the description is not unique. The relationship between the various possible Heegaard splittings that a single manifold may have is often difficult to understand.

Much progress on understanding Heegaard splittings has occurred in the last decade, beginning with the work of Casson and Gordon [CG]. For example, Heegaard splittings of Seifert manifolds can be characterized in a way that makes the classification problem at least accessible (see [MS]). Something of the complexity of Heegaard splittings of hyperbolic manifolds is demonstrated in [LM]. For manifolds which are split up by essential tori into Seifert and hyperbolic pieces, there is some understanding of the relationship between the pieces that can occur and the Heegaard splitting (see [SS]). The remaining (very special) geometric structure which a closed 3-manifold may possess is that of a solvmanifold. Here we completely characterize, up to isotopy, the possible irreducible Heegaard structures on orientable 3-manifolds of this last type.

We will show that, except for precisely two solvmanifolds (those whose monodromy has trace $\pm 3$), any two irreducible Heegaard splittings of the same solvmanifold are isotopic. Depending on the type of solvmanifold, the genus of the Heegaard splitting

Research supported in part by a National Science Foundation grant.

will be either 2 or 3. In each of the two exceptional cases there are exactly two non-isotopic splittings.

Some of what appears here, including a characterization of when a solvmanifold has Heegaard genus two and, for this case, uniqueness up to *homeomorphism* (via examination of the quotient of the hyperelliptic involution), was developed earlier by Takahashi and Ochiai (see [TO]) and Sakuma (see [Sa]). Our viewpoint will be somewhat different and takes advantage of the notion, developed by Casson and Gordon [CG], of strong irreducibility. Michel Boileau offered several useful comments, including one pointing us toward the main argument used here in the proof of Theorem 4.2. We are also indebted to Darren Long and Alan Reid for advice on the special difficulties that arise when $m = 3$.

## 2. Review of Heegaard splittings

**Definition 2.1.** *A* compression body *$H$ is a connected 3-manifold obtained from a closed surface $\partial_- H$ by attaching 1-handles to $\partial_- H \times \{1\} \subset \partial_- H \times I$. Dually, a compression body is obtained from a connected surface $\partial_+ H$ by attaching 2-handles to $\partial_+ H \times \{1\} \subset \partial_+ H \times I$ and 3-handles to any 2-spheres thereby created. The cores of the 2-handles are called* meridian disks

A *Heegaard splitting* $M = A \cup_S B$ of a compact orientable 3-manifold consists of an orientable surface $S$ in $M$, together with two compression bodies $A$ and $B$ so that $S = \partial_+ A = \partial_+ B$ and $M = A \cup_S B$. $S$ itself is called the splitting surface. The genus of the splitting is defined to be the genus of $S$.

A *stabilization* of $A \cup_S B$ is the Heegaard splitting obtained by adding to $A$ a regular neighborhood of a proper arc in $B$ which is parallel in $B$ to an arc in $S$. A stabilization has genus one larger and, up to isotopy, is independent of the choice of arc in $B$. If the construction is done symmetrically to an arc in $A$ instead, the two splittings are isotopic.

Recall the following (see e. g. [CG]): If there are meridian disks $D_A$ and $D_B$ in $A$ and $B$ respectively so that $\partial D_A$ and $\partial D_B$ intersect in a single point in $S$, then $A \cup_S B$ can be obtained by stabilizing a lower genus Heegaard splitting. We then say that $A \cup_S B$ is *stabilized*. If there are meridian disks $D_A$ and $D_B$ in $A$ and $B$ respectively so that $\partial D_A$ and $\partial D_B$ are disjoint in $S$, then (see [CG]) $A \cup_S B$ is *weakly reducible*. If there are meridian disks so that $\partial D_A = \partial D_B$, then $A \cup_S B$ is *reducible*. It is easy to see that reducible splittings are weakly reducible and that (except for the genus one splitting of $S^3$) any stabilized splitting is reducible. It is a theorem of Haken [Ha] that any Heegaard splitting of a reducible 3-manifold is reducible and it follows from a theorem of Waldhausen [W] that a reducible splitting of an irreducible manifold is stabilized. It is an important theorem of [CG] that if a splitting is irreducible but weakly reducible, then maximal simultaneous compression into $A$ and $B$ will create an incompressible surface in $M$.

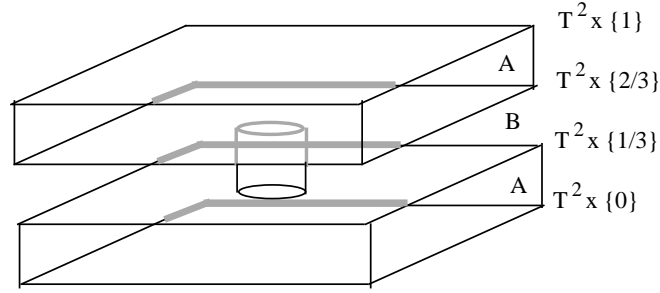

FIGURE 1.

## 3. SOME RELEVANT EXAMPLES

Consider Heegaard splittings of $T^2 \times I$. An easy splitting $T^2 \times I = A \cup_S B$ is obtained by taking $S = T^2 \times \{1/2\}$. This surface divides $T^2 \times I$ into two homeomorphs of $T^2 \times I$, each of which can be regarded as a trivial compression body. In this splitting the components of $\partial T^2 \times I$ lie in different compression bodies.

The easiest splitting for which both components of $\partial(T^2 \times I)$ lie in the same compression body, is given by the following construction: Begin with two copies of the torus, $S' = T^2 \times \{1/3, 2/3\}$ dividing $T^2 \times I$ into three product regions, and set $A' = T^2 \times [0, 1/3] \cup T^2 \times [2/3, 1]$ and $B' = T^2 \times [1/3, 2/3]$. For $E$ a disk in $T^2$, add a vertical tube $E \times [1/3, 2/3]$ to $A'$ (so deleting it from $B'$). This changes $A'$ to a compression body $A$ and, by deleting the tube, changes $B'$ into a genus two handlebody $B$. See Figure 1.

This construction could be generalized to give higher genus Heegaard splittings. For example, begin with the union of three parallel tori $S' = T^2 \times \{1/4, 1/2, 3/4\}$, dividing $T^2 \times I$ into four product regions, and set $A' = (T^2 \times [0, 1/4]) \cup (T^2 \times [1/2, 3/4])$ and $B' = (T^2 \times [1/4, 1/2]) \cup (T^2 \times [3/4, 1])$. Then choose disjoint disks $E_A$ and $E_B$ in $T^2$ and attach the tube $E_A \times [1/4, 1/2]$ to $A'$ and the tube $E_B \times [1/2, 3/4]$ to $B'$ (simultaneously deleting the tubes from $B'$ and $A'$ respectively). The result is a genus three Heegaard splitting $T^2 \times I = A \cup_S B$, but it is easy to see that it is reducible: Let $c_A$ and $c_B$ be simple closed curves in $T^2$ with three properties: They intersect in a single point, they pass through $E_B$ and $E_A$ respectively, and they are disjoint from $E_A$ and $E_B$ respectively. Then $c_B \times [1/4, 1/2]$ and $c_A \times [1/2, 3/4]$ are annuli in $B'$ and $A'$ respectively that intersect in a single point (in $T^2 \times \{1/2\}$). When the tubes are deleted to create $A$ and $B$, the annuli become disks $D_B$ and $D_A$ in $B$ and $A$ respectively, and these disks intersect in a single point. (See Figure 2.) So $T^2 \times I = A \cup_S B$ is a stabilized, hence a reducible splitting.

It is an important and sophisticated theorem of Boileau and Otal ([BO]) that in fact the two elementary examples first given are the only irreducible splittings of $T^2 \times I$. The argument that the third and last example is reducible is particularly important as we now discuss Heegaard splittings of closed orientable 3-manifolds that are torus bundles over the circle.

We review elementary facts and notation about torus bundles over the circle. Any matrix $L \in SL(2, Z)$ induces an orientation preserving homeomorphism on $T^2 =$ $R^2/Z^2$. Unless it is critical, we will not distinguish between $L$ as a matrix and $L$ as a homeomorphism. Points on $T^2$ will be parameterized by ordered pairs $(x_1, x_2), x_i \in R/Z$. Isotopy classes of simple closed curves on $R^2$ are distinguished by their slopes. A circle with slope $n/m$ in $T^2$ is the projection of a line in $R^2$ with slope $n/m$ which in turn is determined by the vectors $\pm \begin{pmatrix} m \\ n \end{pmatrix}$. So in this way we can regard vectors of the form $\begin{pmatrix} m \\ n \end{pmatrix}$, with $m, n$ relatively prime, as parameterizing simple closed curves. With this convention, the minimal number of points in which two curves $c = \begin{pmatrix} m \\ n \end{pmatrix}$ and $c' = \begin{pmatrix} p \\ q \end{pmatrix}$ intersect is given simply by $|det \begin{pmatrix} m & p \\ n & q \end{pmatrix}|$, and this is denoted $c \cdot c'$.

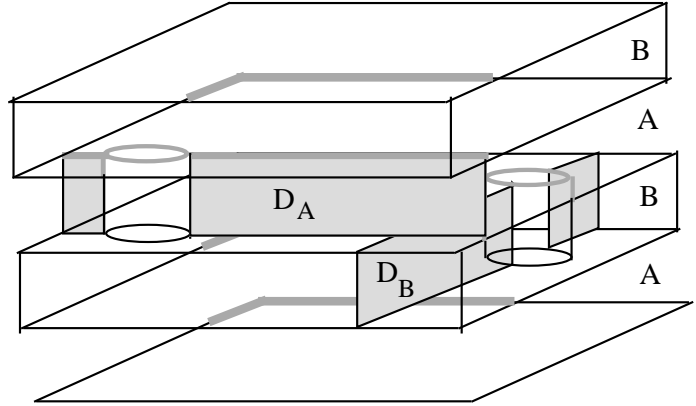

Figure 2.

Let $M_L$ be the mapping cylinder of $L : T^2 \to T^2$. Put another way, $M_L$ is the quotient of $T^2 \times R$ under the identification $(x_1, x_2, t) \sim (L(x_1, x_2), t+1)$ or, alternatively, the quotient of $T^2 \times I$ under the identification $(x_1, x_2, 0) \sim (L(x_1, x_2), 1)$. The homeomorphism $L$ is called the *monodromy* of $M_L$.

The Heegaard splittings of $T^2 \times I$ discussed above suggest Heegaard splittings of $M_L$. Choose disjoint disks $E_A$ and $E_B$ and curves $c_A$ and $c_B$ in the torus $T^2$ as described in the third example above, making sure also that $E_B$ is disjoint from $L(E_A)$. Let $A' \subset M_L$ be the image of $T^2 \times [0, 1/2]$ and $B' \subset M_L$ be the image of $T^2 \times [1/2, 1]$. Attach the tube $E_A \times [1/2, 1]$ to $A'$ and the tube $E_B \times [0, 1/2]$ to $B'$, simultaneously deleting them from $B'$ and $A'$ respectively. The result is a genus three splitting $M_L = A \cup_S B$. Call this the *standard genus three splitting* of $M_L$. The splitting is weakly reducible, since $E_A \times \{3/4\}$ and $E_B \times \{1/4\}$ are disjoint meridian disks.

The disks $D_A \subset A$ and $D_B \subset B$, constructed as in the third example above, will intersect in at least one point (in $T^2 \times \{1/2\}$), but may also intersect in $T^2 \times \{0\} = T^2 \times \{1\}$. Indeed, it will be possible to isotope them so that they are disjoint on $T^2 \times \{0\}$ if and only if $c_B = L(c_A)$. Put succinctly: the standard genus three splitting of $M_L$ is reducible if there is a curve $c_A$ in $T^2$ so that $c_A \cdot L(c_A) = 1$.

We could generalize this construction, as we did that for $T^2 \times I$, by starting with more tori, e. g. the images of the four tori $T^2 \times \{0, 1/4, 1/2, 3/4\}$. But, just as we argued in the third example above, the result is always a reducible splitting.

**Proposition 3.1.** *Suppose $M_L$ is a closed orientable torus bundle over the circle, with mondoromy given by $L \in SL(2, Z)$. The standard genus three splitting is weakly reducible. If the splitting is irreducible, then $c \cdot L(c) \neq 1$ for any simple closed curve $c \subset T^2$.*

$\square$

One way to construct a fiber-preserving homeomorphism from $M_L$ to itself is to begin with a matrix $K \in SL(2, Z)$ that commutes with $L$. Then the automorphism

$$K \times 1_R : T^2 \times R \to T^2 \times R$$

commutes with the covering translation $((x_1, x_2), t) \to (L(x_1, x_2), t+1)$ so it induces a well-defined automorphism of the quotient space $M_L$. We will denote such an automorphism by $\overline{K} : M_L \to M_L$. A useful example later will be the *central* involution $\overline{-I}$ corresponding to the central element $-I$ of $SL(2, Z)$.

The following straightforward observation will be important later:

**Proposition 3.2.** *If $K = L^n$ then $\overline{K} : M_L \to M_L$ is isotopic to the identity.*

**Proof:** The isotopy $h_s : T^2 \times R \to T^2 \times R$ given by $h_s(((x_1, x_2), t)) = ((x_1, x_2), t - ns)$ also commutes with the covering translation and so induces an isotopy $\overline{h_s} : M_L \to M_L$. The maps $h_0$ and $h_1$ are clearly the identity and $\overline{L^n}$ respectively. $\square$

## 4. Solvmanifolds

We now specialize to the case of *solvmanifolds*. These are torus bundles over the circle in which the monodromy $L \in SL(2, Z)$ is Anosov (that is, if it is neither periodic nor does it fix a circle). Equivalent formulations are that $|trace(L)| > 2$ or that $L$ has two irrational eigenvalues or that $L$ is *hyperbolic*. See [Sco].

**Proposition 4.1.** *If $M_L$ is a solvmanifold, then the only irreducible and weakly reducible Heegaard splitting is the standard genus three splitting.*

**Proof:** Suppose $M_L = A \cup_S B$ is a weakly reducible splitting. The main theorem of [CG] shows that, a maximal family of disjoint compressions of $S$ into $A$ and $B$ creates an incompressible surface $T \subset M_L$, and $M_L$ is the union of two distinct submanifolds $A'$ and $B'$ of $M_L$ along $T$. The only incompressible surfaces in $M_L$ are fibers, so $T$ is the union of an even number of fibers and each component of $A'$ and $B'$ is homeomorphic to $torus \times I$.

The reverse construction is easy to describe: $A$ and $B$ can be recovered from $A'$ and $B'$ by attaching tubes through $B'$ and $A'$ respectively. In particular, each component $B_0 \cong T^2 \times I$ of $B'$ has some proper tubes removed in the process of recovering $B$. (Each tube is dual to a meridian disk for $A$ that was compressed in creating $A'$.) The complement of the tubes in $B_0$ then must be a handlebody. In particular, the same tubes could have been removed from $T^2 \times [1/3, 2/3]$ in the second construction of section 3 to give a Heegaard splitting of $T^2 \times I$. But, according to [BO], the result would be reducible unless the removed tubes consist precisely of a single vertical tube. So the result is that in each component of $A'$ (resp. $B'$), a single vertical tube is attached to $B'$ (resp. $A'$) to recover $B$ (resp. $A$). In other words, if $A'$ and $B'$ are each connected, we have the standard genus three splitting; otherwise we have a generalization of this construction, using more fibers, and this was shown in Section 3 to be reducible. □

**Theorem 4.2.** *Suppose a solvmanifold $M_L$ is constructed as the mapping cylinder of $L: T^2 \to T^2$. Then the following are equivalent:*

1. *There is a simple closed curve $c$ in $T^2$ so that $c \cdot L(c) = 1$.*
2. *$L$ is conjugate to a matrix of the form*

$$\begin{pmatrix} \pm m & -1 \\ 1 & 0 \end{pmatrix}, \qquad m \geq 3.$$

3. *Any irreducible splitting of $M_L$ is strongly irreducible.*
4. *Any irreducible splitting of $M_L$ is of genus 2.*
5. *Some irreducible splitting of $M_L$ is of genus 2.*
6. *Some splitting of $M_L$ is strongly irreducible.*

**Proof:** 1) ⇔ 2): If there are curves $c$ and $L(c)$ which intersect in a single point, vectors $v$ and $v'$ in $Z^2$ corresponding to $c$ and $L(c)$ form a basis for $Z^2$ and their signs can be chosen so that $L(v) = -v'$. Then the matrix of $L$ with respect to the basis $(v', v)$ is of the required form. The reverse implication is immediate.

1) ⇒ 3) This follows immediately from 4.1 and 3.1.

4) ⇒ 5) Obvious.

5) ⇒ 6) It is easy to see that a weakly reducible genus two splitting of *any* 3-manifold is reducible.

3) ⇒ 4) and 6) ⇒ 1). We will show that if $M_L$ has a strongly irreducible splitting $A \cup_S B$, then the splitting is of genus two and there is a curve $c$ in $T^2$ so that $c \cdot L(c) = 1$.

The argument is an easy variation of the central argument of [RS], to which we refer for details. Here we present only a sketch. Inside the handlebody $A$ (resp $B$) there is a 1-complex $\Sigma_A$ (resp $\Sigma_B$) to which $A$ (resp $B$) deformation retracts. Then $M_L - (\Sigma_1 \cup \Sigma_2)$ is just a product $S \times \text{int}(I)$. This parameterization of $M_L - (\Sigma_1 \cup \Sigma_2)$ is sometimes called a "sweep-out" by $S$ since $S$ sweeps between one spine and the other.

To each point $(z, t) \in (S^1 \times I)$ we can associate a positioning in $M_L$ of a torus (namely the fiber $T_z$ over $z$) and of $S$ (namely its position $S_t$ during the sweep-out at time $t$). For $t = 0, 1$ we use the spines $S_0 = \Sigma_A$ or $S_1 = \Sigma_B$ respectively instead of a copy of $S$. Put the spines and the sweep-out in general position with respect to the fibering. The result is that at generic points of $(z, t) \in S^1 \times I$ the surfaces $T_z$ and $S_t$ are transverse. Over a 1-dimensional complex $\Gamma$ in $S^1 \times I$ (called the graphic) the surfaces have a single point of tangency, and over the vertices of this complex there are either two points of tangency or a "birth-death" singularity. These last play no important role in the argument.

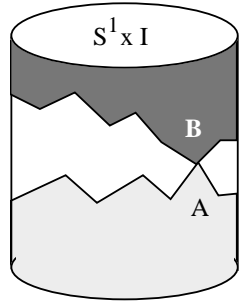

FIGURE 3.

The goal is to understand closed curves of intersection of $S_t$ and $T_z$. We can ignore curves that are inessential in both, and, since $T$ is incompressible, no curve can be essential in $T_z$ and inessential in $S_t$. The important curves will be those that are essential in $S_t$ and inessential in $T_z$. A generic positioning that gives rise to such curves can be labelled $A$ or $B$ depending on whether there is a circle of intersection that bounds a meridian disk of $A$ or $B$. Only one label can occur, since $A \cup_S B$ is strongly irreducible, and, indeed, if two generic regions of $S^1 \times I$ share a common 1-dimensional edge of the graphic $\Gamma$, they can have at most one label $A$ or $B$. (Again, see [RS] for details.)

Next observe what happens when one passes across an edge in $\Gamma$, from one unlabelled region of $(S^1 \times I) - \Gamma$ to another. Since neither $A$ nor $B$ contains an incompressible surface, each $T_z$ must intersect each $S_t$. Because both regions are unlabelled, in each the intersection curves $S_t \cap T_z$ must be essential in $T_z$. But passing from one region to the other through an edge in $\Gamma$ changes the intersection via at most one critical point, and this can't change the slope of the curves $S_t \cap T_z$ (but could only transform a pair of curves into an inessential curve). So, although the number of curves may change, their slope in $T_z$ does not.

This has the following consequence: There cannot be an essential circle $\sigma \subset S^1 \times I$ with the properties that it passes only through unlabelled regions and avoids all vertices of the graphic. For if such a curve existed then, following the remarks of the previous paragraph, one could sweep a torus all the way around the circle (maneuvering $S_t$ as one goes) always keeping the slope of intersection with $S_t$ constant. This would imply that the monodromy would fix a slope on $T^2$, a contradiction.

So now consider what labellings must appear. Near $S^1 \times \{1\}$ each $S_t$ will be very near $\Sigma_A$ and so will intersect any $T_z$ so that some curve of intersection is a meridian of $A$. Similarly, near $S^1 \times \{1\}$ there will always be a curve of intersection bounding a meridian of $B$. In other words, near the former all generic regions will be labelled $A$ and near the latter $B$. On the other hand, we have just seen that two adjacent regions cannot be labelled $A$ and $B$, and there cannot be an essential circle $\sigma \subset S^1 \times I$ with the properties that it passes only through unlabelled regions and avoids all vertices of the graphic. So this means that there must be four regions meeting at a vertex of the graphic, two opposite ones labelled $A$ and $B$ and the other two unlabelled. (See Figure 3.)

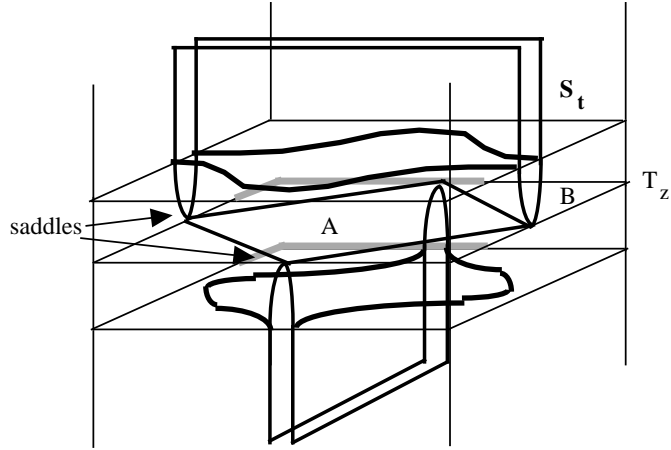

FIGURE 4.

This describes a specific situation (see the remarks preceding [RS, Lemma 5.6]) which, at the vertex of the graphic, can be described as illustrated in Figure 4 ambiently (cf. Figure 5 for the pictures of $S_t$ in $T_z$ in the four different quadrants). $S_t$ and $T_z$ have two points of tangency, connected by four arcs of intersection, to form a diamond in $T_z$. Of the four possible generic intersections that can be created by perturbing the critical points, two give rise to single closed curves, meridians of $A$ and $B$ respectively, and two give rise to pairs of curves with slope $c$ and $c'$ respectively, $c \cdot c' = 1$. Put another way, a small collar neighborhood $\eta(T)$ of $T_z$ in $M_L$ will contain meridian disks of both $A$ and $B$. The bottom of the collar will intersect $S_t$ in curves of slope $c$ and the top in curves of slope $c'$.

Now consider how $S_t$ intersects the closure $M_-$ of $M_L - \eta(T)$, which is also homeomorphic to $T^2 \times I$. $S_t \cap M_-$ cannot be compressible, since any compressing disk would be disjoint from both meridian disks, for $A$ and for $B$, that lie in $\eta(T)$ and this would violate strong irreducibility. Hence $S_t \cap M_-$ cannot be $\partial$-compressible in $M_-$, since $\partial M_-$ consists of tori. (We can ignore or remove $\partial$-parallel annuli in $S_t \cap M_-$.) Hence $S_t \cap M_-$ consists entirely of spanning annuli. The only way these spanning annuli can attach to the appropriate curves in $S_t \cap \eta(T)$ is if $L$ carries $c$ to $c'$, as required. □

It is a formal consequence of Theorem 4.2 that, when condition 1) applies, $M_L$ has a genus two splitting. This can be demonstrated explicitly (see also [TO, Proposition 3]): Let $\Sigma_A \subset M_L$ be the join of two circles: $c \times \{0\} \subset T^2 \times \{0\}$ and the quotient in $M_L$ of the vertical line $\{0,0\} \times R$. Then a neighborhood $A$ of $\Sigma_A$ is clearly a genus two handlebody in $M_L$.

Less obvious is the fact that $B = M_L - \text{int}(A)$ is also a genus two handlebody. To see this, note that $B$ can be obtained as follows: Remove a neighborhood of a vertical spanning arc from $T^2 \times I$. The result is a genus two handlebody $H$. Now glue an annular neighborhood of a curve $c \subset (T^2 \times \{0\} \cap \partial H)$ (disjoint from the neighborhood of the vertical spanning arc that has been removed) to an annular neighborhood of a similar curve $c' \subset (T^2 \times \{1\} \cap \partial H)$ whose slope intersects that of $c$ in a single point.

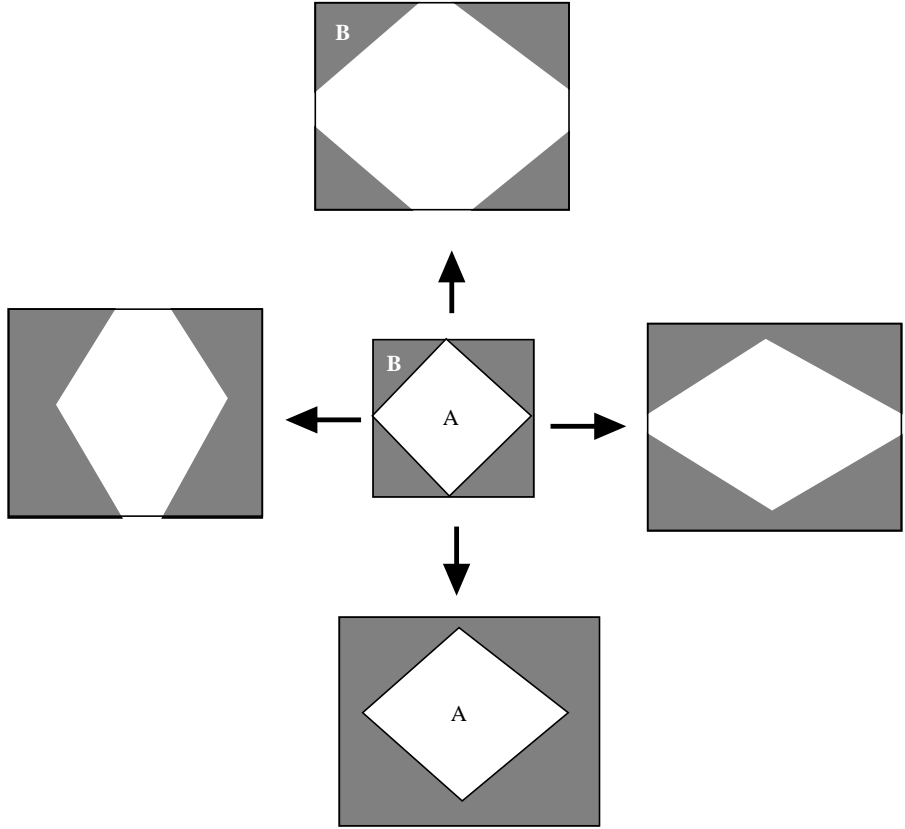

Figure 5.

This last property makes it easy to find disjoint meridian disks $D$ and $D'$ so that $D$ (resp. $D'$) in $H$ intersects $c$ (resp. $c'$) in a single point and is disjoint from $c'$ (resp. $c$). Hence $B$ is obtained from the genus two handlebody $H$ by identifying a longitude of one solid torus summand of $H$ with a longitude of the other solid torus summand. The result is a genus two handlebody.

## 5. Isotopy uniqueness for $trace \neq \pm 3$.

The following algebraic lemma will be essential in verifying when two strongly irreducible splittings of the same solvmanifold are isotopic.

**Lemma 5.1.** *Suppose $m$ is an integer with $|m| \geq 3$ and $K \in SL(2, \mathbb{Z})$ and the matrix*
$$L = \begin{pmatrix} m & -1 \\ 1 & 0 \end{pmatrix}.$$
*If $L$ and $K$ commute then $K = \pm L^n$ for some integer $n$. If instead $KLK^{-1} = L^{-1}$ then $m = \pm 3$.*

**Proof:** There is a direct numerical proof that exploits the simple structure of the eigenvectors of $L$, but we present a geometric argument that we believe is more conceptual and informative.

The modular group $PSL(2,Z)$ acts by Möbius transformations on the upper half plane, $U$, and these are isometries of the hyperbolic metric $ds/y$. Since $trace(L) = m$ and $|m| > 2$ it follows that $L$ is a hyperbolic isometry with some axis $\ell$. Since $K$ commutes with $L$ it follows that $K$ preserves $\ell$ and its orientation. Thus if $K$ is not $\pm I$ then $K$ is also hyperbolic with axis $\ell$. The quotient of the upper half plane by $PSL(2,Z)$ is the modular space $\mathcal{M}$, and we have a branched covering

$$p : U \longrightarrow \mathcal{M}.$$

Now $\mathcal{M}$ is an orbifold which is topologically a disc with two cone points and one cusp. A fundamental domain for this action is

$$D = \{\, x + iy \in U \mid -1/2 \le x \le 1/2 \quad x^2 + y^2 \ge 1 \,\}$$

shown as the shaded region in Figure 6 [M, Theorem 3.2]. It is bounded by two vertical lines and an arc of a circle. The vertical sides of $D$ are identified under $z \mapsto z + 1$ to create the cusp. For each integer $n$ define $C_n$ to be the semi-circle in the upper half plane of radius 1 centered on the point $n$ on the $x$-axis. The circular side, $D \cap C_0$, of $D$ is folded in half by the map $z \mapsto -1/z$. This creates a cone point of order 2 which is the image, $p(i)$, of $i$ and a cone point of order 3 which is the image, $p(exp(\pi i/3))$, of $exp(\pi i/3)$.

We will now describe $\ell$. The Möbius transformation corresponding to $L$ is

$$\tau(z) = m - 1/z.$$

The fixed points of this are

$$z_\pm = \frac{1}{2}\left(m \pm \sqrt{m^2 - 4}\right).$$

Thus $\ell$ is the semi-circle orthogonal to the $x$-axis with endpoints $z_\pm$. The center of this semi-circle is $m/2$. If $m = \pm 3$ then $\ell$ contains the points $1 + i$ and $2 + i$ labelled $a$ and $b$ in Figure 7. These points are in the orbit of $i$ and therefore project on $\mathcal{M}$ to the cone point of order 2. However if $|m| > 3$ then the projection of $\ell$ to $\mathcal{M}$ is disjoint from the cone points. This latter point can be verified as follows.

There is a subarc, $\alpha$, of $\ell$ which is a fundamental domain for the action of $\tau$ on $\ell$. Thus the hyperbolic length of $\alpha$ is the hyperbolic translation length of $\tau$ which is $2cosh^{-1}(|tr(L)|/2)$. Now $z \mapsto -1/z$ maps $C_0$ to itself, reversing the endpoints. Since $\tau$ is the composition of this with $z \mapsto z + m$ it follows that $\tau(C_0) = C_m$ and this map sends the counterclockwise orientation on $C_0$ to the clockwise orientation on $C_m$. The three semi-circles $C_0, \ell, C_m$ are symmetric about $m/2$. We claim that $\alpha$ may be chosen as the subarc of $\ell$ with endpoints on $C_0$ and $C_m$. This is indicated for $m = 4$ in Figure 6. The justification for this as a choice for $\alpha$ is simply that since $\tau(C_0) = C_m$ it follows that $\tau$ does indeed map the endpoint $C_0 \cap \ell$ of $\alpha$ to the other endpoint $C_m \cap \ell$.

When $m = 4$ one checks that $z_-$ is so close to 0 that the endpoint of $\alpha$ on $C_0$ is in the interior of $D \cap C_0$. This is illustrated for $m = 4$ in Figure 6, which has been drawn to scale. Observe that $|z_-|$ is a decreasing function of $|m|$. It then follows that for $|m| \ge 4$ that the endpoint of $\alpha$ on $C_0$ is in the interior of $D \cap C_0$. The situation for $|m| = 3$ is anomalous in this respect, since as shown in Figure 7 the endpoints of $\alpha$ are both outside $D$.

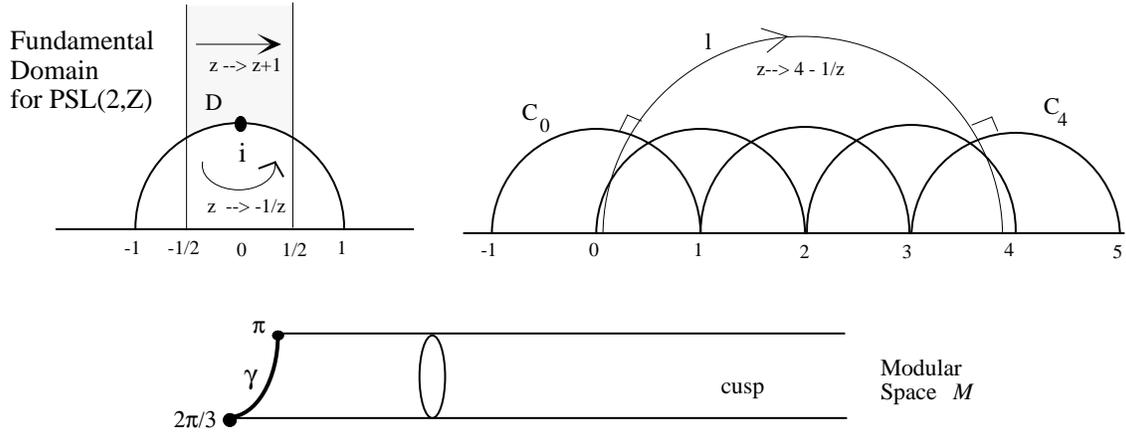

FIGURE 6.

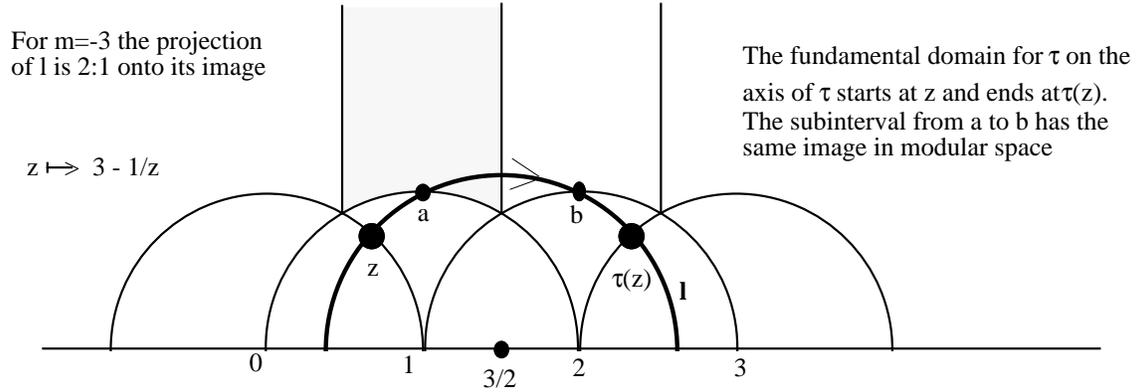

FIGURE 7.

From the symmetry about $m/2$ mentioned before, and the fact that $\tau(C_0) = C_m$ it follows that $C_0, \ell, C_m$ are orthogonal. Hence when $\alpha$ is projected into $\mathcal{M}$ one obtains a geodesic. The fact that $\ell$ is orthogonal to both $C_0$ and $C_m$ ensures that the two ends of $p(\alpha)$ meet on $\mathcal{M}$ along a geodesic arc. (In fact the projection of $\ell$ to $\mathcal{M}$ **must** be a geodesic on $\mathcal{M}$ since the image of a hyperbolic geodesic in a hyperbolic orbifold is always a geodesic, and so the previous reasoning could be omitted.) In the case that $m = \pm 3$ this geodesic is an arc with both endpoints on the cone point of order 2. Otherwise the geodesic is disjoint from both cone points and is thus a closed geodesic.

The orbit of the fundamental domain $D$ under $PSL(2, Z)$ tiles $U$. Define $V$ to be the orbit of $D$ under the group generated by $z \mapsto z + 1$. Then $V$ is the closure of the subset of $U$ lying above the union of the circles $\cup C_n$. The description of the fundamental domain $D$ shows that no element of $PSL(2, Z)$ maps $C_0$ to a vertical side of $D$. Hence the orbit of $C_0$ is disjoint from the interior of $V$. For $|m| \geq 4$ it follows that $\alpha$ meets the orbit of $C_0$ only at its endpoints. The projection $p(C_0)$ of $C_0$ to modular space is the shortest geodesic arc, $\gamma$, on $\mathcal{M}$ connecting the two cone points. Thus the projection of the interior of $\alpha$ is disjoint from $\gamma$. Hence the geodesic, $p(\alpha)$, on $\mathcal{M}$ corresponding to $L$ is **primitive** in other words $p$ is injective on the interior of $\alpha$.

If $K$ and $L$ commute then $K$ also has axis $\ell$ and thus there is an arc $\beta$ of $\ell$ which is a fundamental domain for the action of $K$ on $\ell$. For $|m| \geq 4$ the projection $p(\beta)$ of $\beta$ to $\mathcal{M}$ must wrap around $p(\alpha)$ an integral number of times, say $n$. This implies that the Möbius transformation corresponding to $K$ is $\tau^n$. Since the map from $SL(2, Z)$ to Möbius transformations has kernel $\pm I$ it follows that $K = \pm L^n$. This argument must be modified for $|m| = 3$ since in this case $p$ maps $\alpha$ by a map which is $2 : 1$ onto the arc $p(\alpha)$ in $\mathcal{M}$. The reasoning above then says that $p$ must map $\beta$ an integral number of times over the arc $p(\alpha)$. This means that $length(\beta) = (t/2)length(\alpha)$ for some integer $t$. If $t = 2n + 1$ then $KL^{-n}$ maps $a$ to $b$ and has fundamental domain on $\ell$ the arc with endpoints $a$ and $b$. This is a hyperbolic transformation. It must preserve the tiling of the hyperbolic plane. The shaded translate of $D$ in Figure 7 which contains $a$ must be mapped by this hyperbolic to the adjacent translate of $D$ containing $b$. It is geometrically clear that a hyperbolic element can't do this (but the parabolic $z \mapsto z + 1$ does.) This contradicts that $t$ is odd. Hence $t = 2n$ and for $|m| = 3$ one thus also obtains that $K = L^{\pm n}$.

To prove the final assertion, suppose that $K \in SL(2, Z)$ and that $KLK^{-1} = L^{-1}$. Then $K$ maps $\ell$ to itself reversing endpoints. Hence there is a fixed point, $z$, of $K$ on $\ell$. Thus $K$ is elliptic and $p(z)$ is the cone point of order 2 on $\mathcal{M}$. If $|m| \geq 4$ then $p(\ell)$ is disjoint from the cone points on $\mathcal{M}$ and it follows that $|m| = 3$. □

The proof shows that the case $m = \pm 3$ is different from the other cases because there is a rotation of order 2 in $SL(2, Z)$ with fixed point, $a$, on the axis of $\tau$. This rotation (elliptic) conjugates $\tau$ to its inverse. The Möbius transformation $z \mapsto m - z$ is a half-rotation about $m/2$ which swaps the upper and lower half planes and also conjugates $\tau$ to its inverse. Therefore the composition of these rotations commutes with $\tau$. The square of this composition is $\tau$. This leads to:

**Corollary 5.2.** *Suppose that $L$ is as in the previous lemma and $K \in GL(2, Z)$ and $K$ commutes with $L$. If $det(K) = 1$ then $K = \pm L^n$. If $det(K) = -1$ then $m = \pm 3$ and*

$$m = 3 \implies K = \pm \begin{pmatrix} -2 & 1 \\ -1 & 1 \end{pmatrix} L^n \qquad m = -3 \implies K = \pm \begin{pmatrix} 2 & 1 \\ -1 & -1 \end{pmatrix} L^n.$$

**Proof:** We may assume that $det(K) = -1$, for otherwise the lemma gives the result. The Möbius transformation $z \mapsto m - z$ swaps the endpoints of $\ell$. It corresponds

to a matrix

$$A = \begin{pmatrix} -1 & m \\ 0 & 1 \end{pmatrix}$$

in $GL(2,Z)$ with determinant $-1$. Thus $AK$ is in $SL(2,Z)$ and maps $\ell$ to itself reversing endpoints. Hence $AK$ conjugates $L$ to $L^{-1}$. By the lemma it follows that $m = \pm 3$.

Suppose now that $m = 3$. Let $G$ be the subgroup of $GL(2,Z)$ consisting of all elements which commute with $L$. The intersection of $G$ with $SL(2,Z)$ is a subgroup, $H$, of index 2 and consists of $\pm L^n$. It is readily verified that

$$B = \begin{pmatrix} -2 & 1 \\ -1 & 1 \end{pmatrix}$$

commutes with $L$, and that $det(B) = -1$. Thus $G$ is the union of the two cosets $H$ and $BH$. A similar analysis applies to the case $m = -3$. □

**Theorem 5.3.** *If the monodromy $L$ of a solvmanifold $M_L$ has $|trace(L)| > 3$, then any two irreducible Heegaard splittings of $M_L$ are isotopic.*

**Proof:** We have seen that there are two types of solvmanifolds, those whose monodromy matrix can be written $L = \begin{pmatrix} m & -1 \\ 1 & 0 \end{pmatrix}$ and those that cannot. For those that cannot, it follows from 4.2 that any irreducible splitting is weakly reducible, and from 4.1 that the splitting is then the standard genus three splitting, which is unique up to isotopy.

So we focus entirely on the case in which

$$L = \begin{pmatrix} m & -1 \\ 1 & 0 \end{pmatrix}, \qquad |m| \geq 4.$$

In this case, the discussion following Theorem 4.2 describes a particular genus two splitting, in which $\Sigma_A$ is the join of two circles: $\lambda$, which is the quotient in $M_L$ of the vertical line $\{0,0\} \times R \subset T^2 \times R$, and

$$\begin{pmatrix} 0 \\ 1 \end{pmatrix} \times \{0\} \subset T^2 \times \{0\}.$$

Any other irreducible splitting $M_L = X \cup_Q Y$ must, by Theorem 4.2, be strongly irreducible and, by the last argument in the proof of that theorem, can then be isotoped so that the spine $\Sigma_X$ of one of its handlebodies is the join of $\lambda$ and the curve

$$\begin{pmatrix} a \\ b \end{pmatrix} \times \{0\} \subset T^2 \times \{0\},$$

where the curves $\begin{pmatrix} a \\ b \end{pmatrix}$ and $\begin{pmatrix} u \\ v \end{pmatrix} = L \begin{pmatrix} a \\ b \end{pmatrix}$ intersect in a single point. In particular $\overline{K}(\Sigma_A) = \Sigma_X$.

In other words the matrix

$$K = \begin{pmatrix} u & a \\ v & b \end{pmatrix}$$

is unimodular and conjugates $L$ to a matrix of the form
$$KLK^{-1} = \begin{pmatrix} m' & -1 \\ 1 & 0 \end{pmatrix}.$$
But since trace is preserved by conjugation, $m' = m$ and so $K$ and $L$ commute.

It follows from Cor. 5.2 that $\pm K$ is a power of $L$. Hence by Proposition 3.2, $\overline{\pm K}$ is isotopic to the identity, and so $\Sigma_A$ is isotopic to $\Sigma_X$. □

## 6. When trace is $\pm 3$

The exceptional cases $m = \pm 3$ can be conjugated to
$$\pm \begin{pmatrix} 2 & 1 \\ 1 & 1 \end{pmatrix}.$$
These are well known as the monodromies of the manifolds obtained by zero Dehn-filling the figure eight knot complement and its sister respectively.

In fact these are the only two solvmanifolds with trace $\pm 3$ since there is only one conjugacy class in $SL(2, C)$ with each of these traces. This seems to be well-known (e. g. see [Ra, Section 14]) but we present a geometric proof, shown to us by D. Long and A. Reid, which is in the spirit of the preceding arguments.

**Lemma 6.1.** *There is only one conjugacy class in $SL(2, Z)$ for each of the traces $\pm 3$.*

**Proof:** We consider the case that $trace(A) = 3$, referring again to Figure 7. Let $\ell$ be the axis of $A$, then $p(\ell)$ must intersect $p(C_0)$. Otherwise $p(\ell)$ is contained in the cusp of $\mathcal{M}$ which implies that $A$ is parabolic, contradicting $trace(A) = 3$. This means that we may conjugate $A$ so that $\ell$ intersects $D \cap C_0$. Set
$$A = \begin{pmatrix} a & b \\ c & d \end{pmatrix}.$$
The fixed points of the corresponding Möbius transformation are given by
$$z = \frac{az + b}{cz + d} \quad \text{equivalently} \quad cz^2 + (d - a)z - b = 0.$$
The (Euclidean) distance along the $x$-axis between these two fixed points is
$$\frac{\sqrt{(d-a)^2 + 4bc}}{c} = \frac{\sqrt{(d+a)^2 + 4(bc - ad)}}{c} = \sqrt{5}/c.$$
The lowest point on $D \cap C_0$ is $exp(\pi i/3)$ which has $y$-coordinate $\sqrt{3}/2$. Since $\ell$ is a semi-circle of Euclidian radius $|\sqrt{5}/2c|$, it follows that
$$\frac{\sqrt{5}}{2|c|} \geq \frac{\sqrt{3}}{2}$$
hence $|c| \leq \sqrt{5}/\sqrt{3} < 2$. If $c = 0$ then $A$ is upper triangular hence parabolic, which is a contradiction. Thus $c = \pm 1$. Thus
$$A = \begin{pmatrix} a & \pm(-1 + 3a - a^2) \\ \pm 1 & 3 - a \end{pmatrix}.$$

The fixed points are then
$$\pm \frac{2a - 3 \pm \sqrt{5}}{2}.$$
The semi-circle with these endpoints must intersect $C_0$ thus one of these points has absolute value at most 1. Hence $a = 0, 1, 2, 3$. One now checks that these four matrices are conjugate. The case $trace(A) = -3$ now follows from multiplication by the central element $-I$. □

**Theorem 6.2.** *If the monodromy $L$ of a solvmanifold $M_L$ has $trace(L) = \pm 3$, then $M_L$ has precisely two isotopy classes of irreducible Heegaard splittings. These are strongly irreducible, genus two, and the product of their associated hyperelliptic involutions is the central involution.*

**Proof:** According to Lemma 6.1, we may as well take
$$L = \begin{pmatrix} \pm 3 & -1 \\ 1 & 0 \end{pmatrix},$$
so much of the argument used in proving Theorem 5.3 applies. The cases are symmetric, so we will take $m = 3$, and, following the argument of Theorem 5.3, we need to consider exactly the case when $det(K) = -1$ and $K = \pm \begin{pmatrix} -2 & 1 \\ -1 & 1 \end{pmatrix}$. In particular, a possibly alternative Heegaard splitting is one which replaces the simple closed curve $\begin{pmatrix} 0 \\ 1 \end{pmatrix} \times \{0\} \subset \Sigma_A$ with $\begin{pmatrix} 1 \\ 1 \end{pmatrix} \times \{0\} \subset T^2 \times \{0\}$. We will call the resulting spine of the possibly alternate Heegaard splitting $\Sigma_X$.

In order to examine these two potentially different Heegaard splittings more carefully, it will be convenient to isotope $\Sigma_X$ a bit and introduce some helpful notation. Recall that $\lambda$ is the "vertical" circle in $M_L$ that is the quotient of the line $\{0, 0\} \times R$. Let $\alpha \subset T^2$ be the curve $\begin{pmatrix} 0 \\ 1 \end{pmatrix}$, $\beta \subset T^2$ be the curve $\begin{pmatrix} 1 \\ 1 \end{pmatrix}$ and $\gamma \subset T^2$ be the curve $\begin{pmatrix} 2 \\ 3 \end{pmatrix}$. Then we take, as before, $\Sigma_A = \lambda \cup (\alpha \times \{0\})$ and move $\Sigma_X$ to $\lambda \cup (\beta \times \{1/2\})$.

Let $\rho : T^2 \to T^2$ be the involution $\rho(x_1, x_2) = (x_2, x_1)$. Since $\rho L \rho = L^{-1}$, or, equivalently $(\rho L)^2 = I$, the orientation preserving involution $\hat{\rho}$ of $T^2 \times I$ given by $\hat{\rho}(x_1, x_2, t) = (x_2, x_1, 1 - t)$ descends to an involution of $M_L$. The fixed point set of $\hat{\rho} : M_L \to M_L$ consists of two circles $\beta \times \{1/2\} \subset T^2 \times \{1/2\}$, and $\gamma \times \{0\} \subset T^2 \times \{0\}$. The first is obvious, and the latter follows since $L(\gamma) = \rho(\gamma)$. Similarly, note that $\hat{\rho}$ takes both $\Sigma_X$ and $\Sigma_A$ to themselves, the latter because $\rho(\alpha) = -L(\alpha)$. The first fixed curve $\beta \times \{1/2\} \subset T^2 \times \{1/2\}$ intersects the spine $\Sigma_A$ in a single point and the second $\gamma \times \{0\} \subset T^2 \times \{0\}$ intersects it twice. So $\hat{\rho}$ is the hyperelliptic involution induced by the Heegaard splitting $A \cup_S B$. Note that it preserves the splitting $X \cup_Q Y$.

If we replace $\rho$ by $-\rho$ in the above argument, much remains the same. Again we get an involution $\hat{\rho}' : M_L \to M_L$ but now instead of *fixing* $\beta \times \{1/2\}$, $\hat{\rho}'$ reflects it and so there are two fixed points on that circle. Similarly, instead of reflecting $\alpha \times \{1/2\}$, $\hat{\rho}'$ fixes it. In other words $\hat{\rho}'(x_1, x_2, t)$ is the hyperelliptic involution on $\Sigma_X$ and fixes $\Sigma_A$. The product $\hat{\rho} \cdot \hat{\rho}' = -\overline{I}$ is the central involution.

All that remains is to show that the splittings $A \cup_S B$ and $X \cup_Q Y$ described above are not isotopic. But clearly isotopic Heegaard splittings will have isotopic hyperelliptic involutions, so their product will be isotopic to the identity. But the product of these hyperelliptic involutions is

$$-\overline{I} : M_L \to M_L$$

and this cannot be isotopic to the identity. For if it were, the lift of the isotopy to $T^2 \times R$ would force $-I = L^n$ for some $n$, and this is clearly impossible (since an eigenvalue of $L$ is $> 1$). □

## 7. Commensurability Relations between Solvmanifolds.

For three-dimensional solvmanifolds we have seen that if the monodromy has a special form then the Heegaard genus is 2 and otherwise it is 3. In the process we have noted that there is only one solvmanifold of trace 3. In this section we show that such uniqueness is "virtually" true for all other traces. That is, up to taking finite covers, a three-dimensional solvmanifold is determined by the trace of the monodromy.

We will say that two self-homeomorphisms $\phi_1, \phi_2$ of a torus are **virtually conjugate** if there is a self-homeomorphism $\phi$ of a torus which covers both $\phi_1$ and $\phi_2$ under suitable coverings. An equivalent way to say this is that, taking lifts of $\phi_1, \phi_2$ to suitable finite coverings, the lifts are conjugate. This is an equivalence relation.

**Theorem 7.1.** *Suppose that $A, B \in SL(2, Z)$ both have trace equal to $m$ and that $|m| \geq 3$. Then the corresponding linear self-homeomorphisms of the torus are virtually conjugate. Conversely, if two linear self-homeomorphisms of the torus are virtually conjugate, then they have equal traces.*

**Proof:** Suppose that $f, g$ are self-homeomorphisms of the torus corresponding to $A, B$ in $SL(2, Z)$ and that there is a finite cover of the torus for which $g$ covers $f$. Then $f, g$ are both covered by the same linear automorphism of $R^2$. The trace of this linear map equals the trace of both $f$ and $g$.

Conversely, suppose $trace(A) = trace(B) = m$ and that $|m| \geq 3$. Then $A$ and $B$ represent translations in the hyperbolic plane by equal distances. There is an isometry of the hyperbolic plane taking the oriented axis of $A$ to that of $B$. Thus $A$ and $B$ are conjugate in $SL(2, R)$. Thus there is $P$ in $SL(2, R)$ with $PA = BP$. This may be thought of as a linear system of equations with integer coefficients and the entries in $P$ are the unknowns. It follows that there is $P \in GL(2, Q)$ satisfying this equation. We may thus take $P$ to have integer entries and $det(P) \neq 0$. Regarding $P$ as a monomorphism of the group $Z^2$ we see that $\Lambda$, the image of $P$, is a subgroup of finite index in $Z^2$. Since $PA = BP$ it follows that $B(\Lambda) = \Lambda$. Regard $B$ as a homeomorphism of a torus $T$. Let $\tilde{T}$ be the covering of $T$ corresponding to the subgroup $\Lambda$ of $\pi_1 T$. Then $B$ is covered by an automorphism $\tilde{B}$ of $\tilde{T}$. Then $P$ may be regarded as a homeomorphism from $T$ to $\tilde{T}$ which conjugates $A$ to $\tilde{B}$. □

The following provides a strong form of commensurability:

**Theorem 7.2.** *Suppose that $A, B \in SL(2, Z)$ satisfy $|trace(A)|, |trace(B)| \geq 3$. Let $M_A$ and $M_B$ be the solvmanifolds with these monodromies. Then $M_A$ has a finite cover homeomorphic to $M_B$ and vice versa if and only if $trace(A) = trace(B)$.*

**Proof:** If $trace(A) = trace(B)$ then there is a finite cover $\tilde{M}_B$ obtained by taking the finite cover, $\tilde{T}_B$, of the fiber torus, $T_B$, in $M_B$ used in the previous proof. Since $B$ is covered by the map $\tilde{B}$ of $\tilde{T}_B$, this provides the covering $\tilde{M}_B$. It is now clear that $M_A$ is homeomorphic to $\tilde{M}_B$. For the converse, one shows that if $M_A$ covers $M_B$ then there is $n \geq 1$ with $trace(A^n) = trace(B)$. □

Daryl Cooper, Mathematics Department, University of California, Santa Barbara, CA USA

*E-mail address*: `cooper@math.ucsb.edu`

Martin Scharlemann, Mathematics Department, University of California, Santa Barbara, CA USA

*E-mail address*: `mgscharl@math.ucsb.edu`